\documentclass[11pt]{article}
\usepackage{amsmath,amssymb,amsthm,amsfonts,verbatim}
\usepackage{graphicx}
\usepackage{microtype}
\usepackage{url}

\title{Infinite generation of the kernels of the Magnus and Burau
  representations}

\author{Thomas Church and Benson Farb
\thanks{The second author gratefully acknowledges support from the National
Science Foundation.}}

\theoremstyle{plain}
\newtheorem{theorem}{Theorem}[section]
\newtheorem{proposition}[theorem]{Proposition}
\newtheorem{lemma}[theorem]{Lemma}
\theoremstyle{definition}
\newtheorem{corollary}[theorem]{Corollary}

\newcommand{\nc}{\newcommand}
\nc{\dmo}{\DeclareMathOperator}

\nc{\Z}{\mathbb{Z}}
\nc{\Tor}{\mathcal{I}}
\nc{\bwedge}{\textstyle{\bigwedge}}

\dmo\im{im}
\dmo\Sp{Sp}
\dmo\Mod{Mod}
\dmo\CAT{CAT}
\dmo\ab{ab}
\dmo\Hom{Hom}
\dmo\GL{GL}
\dmo\Aut{Aut}
\dmo\Mag{Mag}
\dmo\Bur{Bur}
\dmo\kernel{kernel}

\nc{\figurelantern}{1}
\nc{\figurearcs}{2}
\nc{\figuresurface}{3}
\nc{\figurefourcurves}{4}
\nc{\figurebasis}{5}
\nc{\figuregenustwo}{6}
\nc{\figurebutterfly}{7}

\nc{\coloneq}{\mathrel{\mathop:}\mkern-1.2mu=}
\nc{\co}{\colon\thinspace}

\begin{document}
\maketitle
\begin{abstract}
  Consider the kernel $\Mag_g$ of the Magnus representation of the
  Torelli group and the kernel $\Bur_n$ of the Burau representation of
  the braid group.  We prove that for $g\geq 2$ and for $n\geq 6$ the
  groups $\Mag_g$ and $\Bur_n$ have infinite rank first homology.  As
  a consequence we conclude that neither group has any finite generating
  set.  The method of proof in each case consists of producing a kind
  of ``Johnson-type'' homomorphism to an infinite rank abelian group,
  and proving the image has infinite rank.  For the case of
  $\Bur_n$, we do this with the assistance of a computer calculation.
\end{abstract}

\maketitle
\section{Introduction}
\label{section:introduction}

\noindent
\textbf{The Magnus kernel.}
Let $S\coloneq S_{g,1}$ be a compact, connected, oriented surface of
genus $g\geq 2$ with one boundary component.  Let $\Mod_{g,1}$ denote
the \emph{mapping class group} of $S$, which is the group of homotopy
classes of orientation-preserving homeomorphisms of $S$ which fix
$\partial S$ pointwise.  Let $\Tor_{g,1}$ denote the \emph{Torelli
  group}, which is the subgroup of $\Mod_{g,1}$ consisting of elements
that act trivially on $H\coloneq H_1(S,\Z)$.

$\Mod_{g,1}$ acts on the fundamental group $\pi_1(S)$, inducing an
action on the solvable quotient $\Gamma/\Gamma^3$, where
$\Gamma\coloneq \pi_1(S)$, $\Gamma^2=[\Gamma,\Gamma]$ and $\Gamma^3=
[\Gamma^2,\Gamma^2]$ are the first three terms of the derived series
of $\Gamma$. In this paper we consider the group \[\Mag_g\coloneq
\kernel(\Mod(S)\to \Aut(\Gamma/\Gamma^3)).\] In 1939, Magnus
(\cite{Ma}; see also \cite[Chapter 3]{Bi}) used the Fox calculus to
construct a representation
\[r\co \Tor_{g,1}\to \GL_{2g}(\Z H)\] now called the \emph{Magnus
  representation}. It follows from \cite[Theorem 4.9]{Fox} that the
kernel of $r$ coincides with $\Mag_g$.  This group is called the
\emph{Magnus kernel}.

It was an open question for some time whether or not $\Mag_g$ is
nontrivial.  This was settled in the affirmative by Suzuki in
\cite{Su}. The first main result of this paper is that $\Mag_g$ is in fact
quite large.

\begin{theorem}
\label{theorem:main}
For $g\geq 2$ the group $H_1(\Mag_g,\Z)$ has infinite rank.
\end{theorem}

As the abelianization of a finitely-generated group has finite rank,
we deduce the following.
\begin{corollary}
\label{corollary:main}
For $g\geq 2$ the group $\Mag_g$ has no finite generating set.
\end{corollary}

The idea of our proof of Theorem~\ref{theorem:main} is to define
a kind of ``Johnson-type'' homomorphism (see \cite{J1}):
\[\Psi\co \Mag_g\to \Hom\big(G^{\ab},\bwedge^2 G^{\ab}\big)\]
where $G=[\Gamma,\Gamma]$ and $G^{\ab}$ denotes the abelianization of
$G$.  We then construct infinitely many linearly independent elements
contained in the image.

\bigskip
\noindent
\textbf{The Burau kernel.} Let $B_n$ denote the braid group on
$n$ strands.  $B_n$ can be realized (see Section~\ref{section:braid}
below) as a subgroup of the automorphism group $\Aut(F_n)$ of the free
group of rank $n$.  The \emph{Burau representation} is a homomorphism
\[\rho_n\co B_n\to \GL_{n}(\Z[t,t^{-1}]).\]
We define the \emph{Burau kernel}, denoted $\Bur_n$, to be the kernel
of $\rho_n$.  Let $K$ be the kernel of the homomorphism $F_n\to \Z$
taking each fixed generator of $F_n$ to $1$.  It follows easily from
\cite{Fox} that
\[\Bur_n=\kernel(B_n\to \Aut(F_n/[K,K])).\]

While $\rho_3$ is faithful, it was a longstanding problem as to
whether or not $\rho_n$ is faithful (i.e. whether $\Bur_n$ is nontrivial) for
$n>3$.  This was solved by Moody \cite{Mo}, Long--Paton \cite{LP}, and
Bigelow \cite{Big} in various cases, with the result that $\Bur_n$ is
nontrivial for $n\geq 5$; the case of $n=4$ is still open.  Our next
main result is that $\Bur_n$ is in fact quite large for $n\geq 6$.

\begin{theorem}
  \label{theorem:main2}
  For $n\geq 6$ the group $H_1(\Bur_n,\Z)$ has infinite rank; in
  particular, $\Bur_n$ has no finite generating set.
\end{theorem}

To prove Theorem~\ref{theorem:main2} we construct, similarly to the
proof of Theorem~\ref{theorem:main} above, a homomorphism
\[\Phi\co \Bur_n\to \Hom\big(K^{\ab}, \bwedge^2 K^{\ab}\big).\]

The elements which have been constructed in the kernel of the Burau
representation are geometrically elegant, but algebraically very
complicated; for example, the element of $\Bur_7$ found by Long--Paton
can be described by a single diagram, but as a free group automorphism
sends generators of $F_7$ to words of length up to $475137$. Thus we
need the assistance of a computer in order to calculate $\Phi$
explicitly (see Section~\ref{section:braid} below for a full
discussion). For the computations in this paper we use a simpler
element $\phi_B\in \Bur_n$ for $n\geq 6$ found by Bigelow, which takes
generators to words of length no more than 9841. Once we compute the
form of $\Phi(\phi_B)$, we then use an equivariance property of $\Phi$
to prove that the image of $\Phi$ has infinite rank, from which
Theorem~\ref{theorem:main2} follows.

We remark that, as Problem 6.24 of \cite{Mor}, Morita posed the
problem of determining the kernel of the Magnus and Burau (among
other) representations.  Theorem \ref{theorem:main} and Theorem
\ref{theorem:main2} can be viewed as a partial answer to this problem.
\pagebreak

\noindent
\textbf{Acknowledgements.} We are grateful to William Goldman, whose
Mathematica notebook FreeGroupAutos.nb was very helpful in our
computations of the expression in Appendix A. We would also like to
thank Dan Margalit and Tam Nguyen Phan for careful comments on an
earlier version of this paper.

\section{Defining the homomorphisms}
\label{section:definition}

The following construction works whenever one considers a group of
automorphisms of the universal $2$-step nilpotent quotient of  a group $G$
acting trivially on its abelianization. Johnson \cite{J1} considered
the case $G=\Gamma=\pi_1(S)$.

With $\Gamma$ equal to $\pi_1(S)$ or $F_n$ as in the introduction, we take
$G\coloneq[\Gamma,\Gamma]$ or $G\coloneq K$ respectively.  In either
case, let $G_i$ be the lower central series of $G$, defined
inductively by $G_1=G$ and $G_{i+1}=[G,G_i]$.  Consider the exact
sequence
\begin{equation}
  \label{eq:basic1}
  1\to G_2\to G\to G^{\ab}\to 1.
\end{equation}
Centralizing (\ref{eq:basic1}) gives
\begin{equation}
  \label{eq:key}
  1\to G_2/G_3\to G/G_3\to G^{\ab}\to 1.
\end{equation}
Since $G$ is free, taking (\ref{eq:basic1}) as a presentation for
$G^{\ab}$, Hopf's formula gives that
\[G_2/G_3\approx \bwedge^2 G^{\ab}.\] 

$\Aut(\Gamma)$ acts on $\Gamma$, and thus on $G$, and the isomorphism
$\nu\co G_2/G_3\approx \bwedge^2 G^{\ab}$ respects the action of
$\Aut(\Gamma)$ on both sides. In particular, conjugation by $\Gamma$
descends to an action on $G^{\ab}$ by $H=\Gamma/[\Gamma,\Gamma]$ or by
$\Z=\Gamma/K$ respectively. In the case $G=[\Gamma,\Gamma]$, the fact
that $\Mag_g$ acts trivially on $\Gamma/\Gamma^3$ implies that
$\Mag_g$ acts trivially on $G^{\ab}=\Gamma^2/\Gamma^3$ and on
$\bwedge^2 G^{\ab}$.  Similarly, in the case $G=K$, we have that
$\Bur_n$ acts trivially on $G^{\ab}$ and on $\bwedge^2 G^{\ab}$.

Let $f\in\Mag_g$ (resp.\ $f\in\Bur_n$) be given.  For $x\in G^{\ab}$, pick
any lift $\tilde{x}\in G$.  Since $f$ acts trivially on both the
quotient and kernel of (\ref{eq:key}), we see that
$f(\tilde{x})\tilde{x}^{-1}$ lies in the kernel $G_2/G_3$, which we
identify with $\bwedge^2 G^{\ab}$ via the isomorphism above.  One checks, exactly as in \cite{J1}, that
\[\delta_f\co G^{\ab}\to \bwedge^2 G^{\ab}\]
defined by $\delta_f(x)\coloneq f(\tilde{x})\tilde{x}^{-1}$ is a
well-defined homomorphism; in fact, the resulting map $\delta_f$ is
$\Z H$--linear (resp.\ $\Z[t,t^{-1}]$--linear) with respect to the
conjugation action on $G^{\ab}$. This is equivalent to the claim
that \[\delta_f(\gamma x\gamma^{-1})\equiv \gamma
\delta_f(x)\gamma^{-1}\bmod{G_3},\] which can be checked as
follows. The difference between the left and right side is
\[\big(f(\gamma x\gamma^{-1})\gamma
x^{-1}\gamma^{-1}\big)\big(\gamma
f(x)x^{-1}\gamma^{-1}\big)^{-1}=f(\gamma)f(x)f(\gamma)^{-1}\gamma
f(x)^{-1}\gamma^{-1},\] which is conjugate to $[\gamma^{-1}
f(\gamma),f(x)]$. The condition on $f$ implies that $ f(\gamma)\equiv
\gamma\bmod{G_2}$, so $\gamma^{-1} f(\gamma)\in G_2$ and $[\gamma^{-1}
f(\gamma),f(x)]\in G_3$ as desired.

One also checks, exactly as in \cite{J1}, that in the case
$G=[\Gamma,\Gamma]$, defining the map $\Psi$ by $\Psi(f)\coloneq \delta_f$ gives a
well-defined homomorphism;
\begin{equation}
  \label{eq:johnson1}
  \Psi\co \Mag_g\to \Hom\big(G^{\ab},\bwedge^2 G^{\ab}\big).
\end{equation}
and, in the case $G=K$, defining $\Phi(g)\coloneq \delta_f$ gives a
well-defined homomorphism:
\begin{equation}
  \label{eq:johnson2}
  \Phi\co \Bur_n\to \Hom\big(G^{\ab},\bwedge^2 G^{\ab}\big).
\end{equation} 
The homomorphisms $\Psi$ and $\Phi$ are equivariant with
respect to the natural $\Aut(\Gamma)$-actions on the source and
target.

\section{Computing the image of \boldmath$\Psi$}
\label{section:Psi}
Let $S_{0,4}$ denote the $2$-sphere with $4$ open disks removed.  A
\emph{lantern} in $S$ is an embedding $S_{0,4}\hookrightarrow S$.
Consider the two simple closed curves $\alpha$ and $\beta$ and the
three arcs $A_1,A_2$ and $A_3$ on $S_{0,4}$ given in
Figure~\figurelantern.
\begin{figure}[h]
  \label{figure:curves}
  \begin{center}
    \includegraphics[width=105mm]{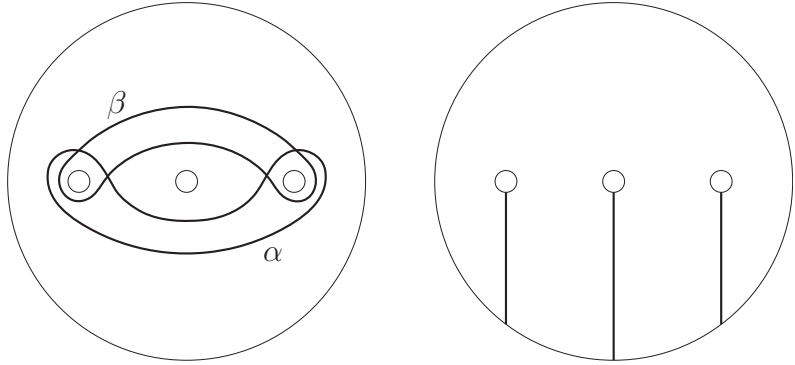}\\
  \end{center}
  \caption{The simple closed curves $\alpha$ and $\beta$, and the
    arcs $A_1,A_2,A_3$.}
\end{figure}

One directly computes the action of $f\coloneq
T_{\alpha}T_{\beta}^{-1}$ on $A_1$, $A_2$ and $A_3$, as follows (see
Figure~\figurearcs).  Let $x$, $y$, and $z$ be the loops which
begin with $A_1$, $A_2$ and $A_3$, respectively, go clockwise around
the appropriate boundary component of $S_{0,4}$, then come back along
the same arc $A_i$.  Let $X,Y,Z$ be the inverses of $x,y,z$ in
$\pi_1(S_{0,4})$.  Then:
\[\begin{array}{l}
  f(A_1)= xyXzxYXZA_1= [xyX,z]A_1\\
  \\
  f(A_2)= ZXzxA_2= [Z,X]A_2\\
  \\
  f(A_3)= ZXzxYXZxzxyXA_3 = [ZXz,xYX]A_3
\end{array}\]

Let $L$ be an embedding of a lantern in $S$ with the property that
each of the four boundary curves of $L$ are separating in
$S$.\footnote{To formally identify $x,y,z$ with elements of
  $\Gamma=\pi_1(S)$, we choose a basepoint on $\partial S$, and arcs
  from this basepoint to $L$ meeting $L$ in one point.  Since $f$ is
  the identity off of $L$, any ambiguity in the choice of these paths
  to $L$ does not affect the computation.}  In this case we can
observe that $T_{\alpha}T_{\beta}^{-1}\in \Mag_g$, as follows. Note
that the elements corresponding to $x,y,z$ all lie in
$\Gamma^2$. Furthermore, $\Gamma=\pi_1(S)$ has a basis where each
element $c$ is either disjoint from $L$, or else of the form
$c=A\gamma A^{-1}$, where $A$ is an arc intersecting $L$ in some $A_i$
and $\gamma$ is a loop disjoint from $L$. In the former case the
element $f=T_{\alpha}T_{\beta}^{-1}$ fixes $c$. In the latter case,
assume for example that $A$ intersects $L$ in $A_2$; then we have
\[f(c)=f(A\gamma A^{-1})=f(A)\gamma f(A)^{-1}=[Z,X]A\gamma
A^{-1}[X,Z]=[Z,X]c[X,Z]\] Since $x,y,z\in \Gamma^2$, we have $[Z,X]\in
\Gamma^3$; thus $f(c)\equiv c\bmod{\Gamma^3}$. The same is true for
$A_1$ and $A_3$, so we conclude that $f(c)\equiv c\bmod{\Gamma^3}$ for
all elements of a basis for $\Gamma$, implying
$T_{\alpha}T_{\beta}^{-1}\in \Mag_g$. Suzuki gave a more illuminating proof
that elements of this form lie in $\Mag_g$ in \cite{Su2}.
\begin{figure}
\label{figure:arcs}
\begin{center}
    \includegraphics[width=162mm]{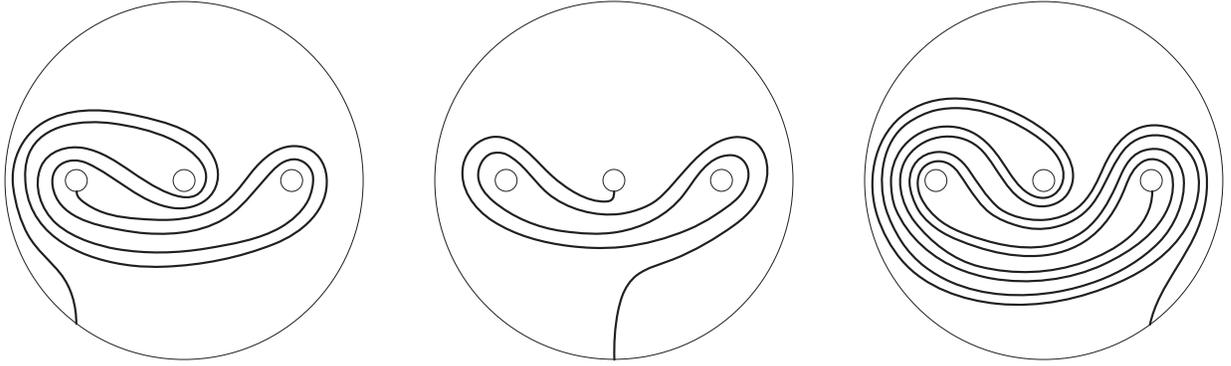}\end{center}
  \caption{The arcs $f(A_1), f(A_2)$ and $f(A_3)$.}
\end{figure}

We are now ready to compute $\Psi$.  For $a,b\in\Gamma$, we denote by
$\{a,b\}$ the image of $[a,b]\in G$ in $G^{\ab}$ under the
abelianization map.

\begin{proposition}\label{prop:lantern} Let $L$ be a lantern embedded
  in $S$ so that each of the four boundary curves of
  $L$ are separating in $S$.  Let $a$ and $b$ be
  loops intersecting $L$ in $A_1$ and $A_2$. Then
\begin{equation}
\label{eq:lantern5}
\Psi(T_{\alpha}T_{\beta}^{-1})(\{a,b\})=(a-1)(b-1)\big[ x\wedge z
  + y\wedge z\big]
  \end{equation}
\end{proposition}

Note that the right hand side of (\ref{eq:lantern5}) is an element of
$\bwedge^2 G^{\ab}$, considered as a $\Z H$-module, and $a,b$ are
taken to be elements of $H$.

\begin{proof}
  As in the computation above, we have
  \[f([a,b])=[f(a),f(b)]=[wa,vb]\] where \[w= [[xyX,z],a] \ \
  \mbox{and}\ \ v=[[Z,X],b].\] From the assumption on the embedding of
  $L$ we have $x,y,z\in G$, and thus $w,v\in G_2$. We will use the
  following commutator identities, which hold in any group; we write
  ${}^xy$ for $xyx^{-1}$.
  \[[wa,b]={}^w[a,b]\ [w,b]
  \qquad\qquad[a,vb]=[a,v]\ {}^v[a,b]\]
  \begin{figure}[h]
    \label{figure:surface}
    \begin{center}
      \includegraphics[width=105mm]{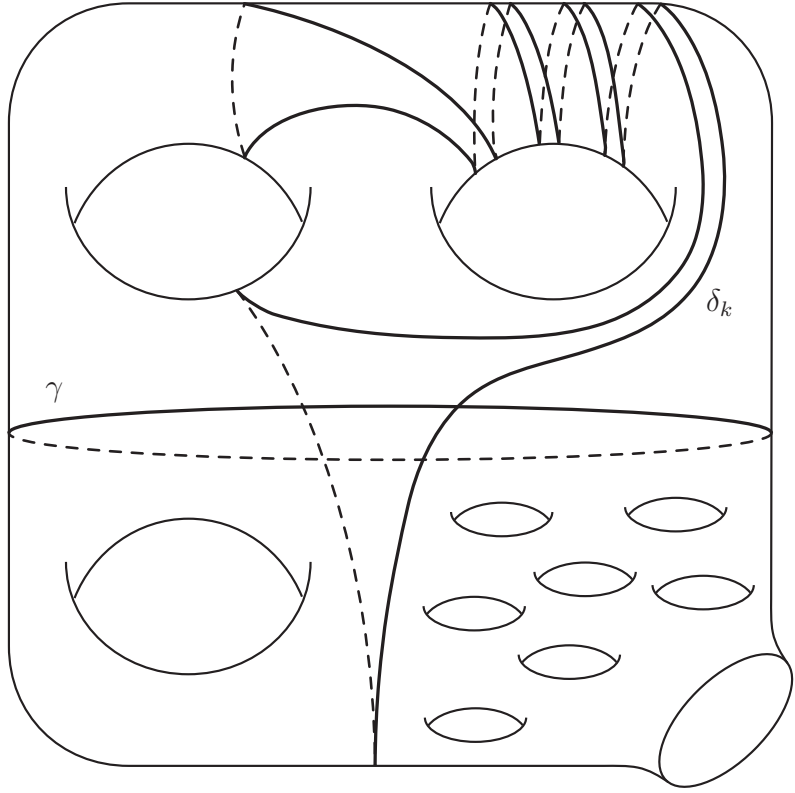}\\
    \end{center}
    \caption{The curves $\gamma$ and $\delta_k$ for $k=3$.}
  \end{figure}
  We then find that
  \[[wa,vb]={}^w[a,v]\ {}^{wv}[a,b]\ [w,v]\ {}^v[w,b]\] 
  Note that the
  second term lies in $G$, the first and fourth in $G_2$, and the
  third in $G_3$.
  
  We want to compute $f([a,b])[a,b]^{-1}$ as an element of
  $G_2/G_3$. Note that ${[w,v]\equiv 0\bmod G_3}$, and that conjugating
  an element of $G$ by an element of $G_2$ is a trivial operation
  modulo $G_3$. Finally, since $[[a,b],[w,b]]\in G_3$, we can move
  $[a,b]$ to the right to cancel $[a,b]^{-1}$. We thus obtain
  \begin{align*}
    f([a,b])[a,b]^{-1}&={}^w[a,v]\ {}^{wv}[a,b]\ [w,v]\ {}^v[w,b]\ [a,b]^{-1}\\
    &\equiv [a,v][a,b][w,b][a,b]^{-1}\bmod G_3\\
    &\equiv [a,v][w,b]\bmod G_3.
  \end{align*} Recall that action of $\Gamma$ on $\Gamma$ by
  conjugation decends to a $\Z H$ action on $G^{\ab}$.  Recall from
  above the isomorphism $\nu\co G_2/G_3 \to \bwedge^2 G^{\ab}$.  Since
  the homology class of $x$ is trivial in $H$, we have
  \[\nu([xyX,z])=y\wedge z \ \ \mbox{and}\ \ \nu([Z,X])=z\wedge x.\]
  It follows that \[\nu(w)=\nu([[xyX,z],a])=(1-a)y\wedge z\]
  and \[\nu(v)=\nu([[Z,X],b])=(1-b)z\wedge x.\]
  We therefore have that
  \[\nu([a,v][w,b])=(a-1)v-(b-1)w=(a-1)(1-b)z\wedge
  x-(b-1)(1-a)y\wedge z.\]
  We conclude that
  \[\Psi(T_{\alpha}T_{\beta}^{-1})(\{a,b\})=(a-1)(b-1)\big[ x\wedge z
  + y\wedge z\big]\]
  as desired.
\end{proof}
\begin{figure}
    \label{figure:fourcurves}
    \begin{center}
      \includegraphics[width=90mm]{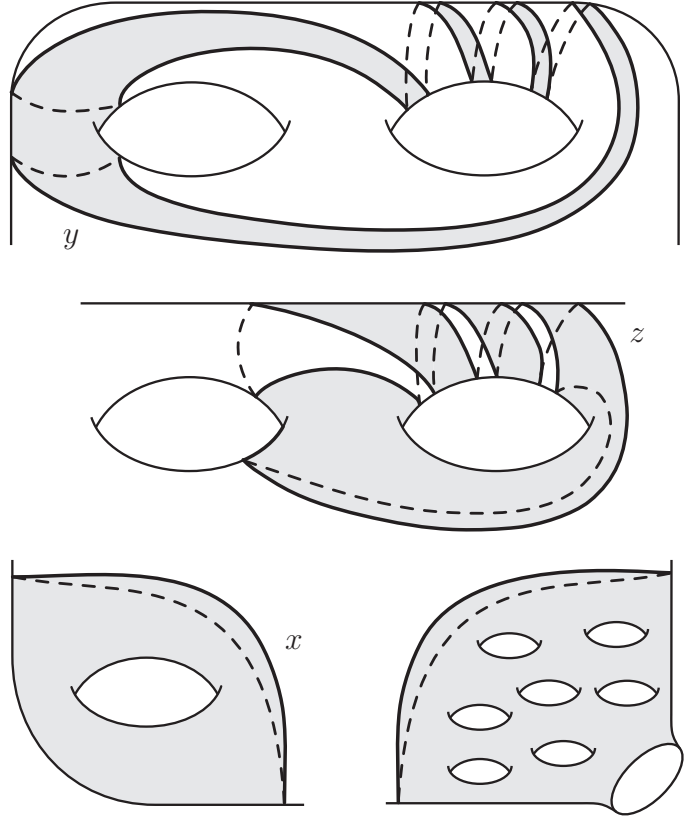}\\
    \end{center}
    \caption{The boundary curves of $L_k$; the subsurfaces cut off by
      these curves are shaded.}
  \end{figure}

\begin{theorem}The image of $\Psi$ has infinite rank for $g\geq 3$.
  \label{theorem:imageA}
\end{theorem}
\begin{proof}
  Let $\gamma$ and $\delta_k$ be the curves depicted in
  Figure~\figuresurface. The figure depicts the case $k=3$; in general
  $\delta_k$ has $k$ twists around the upper right
  handle. (Specifically, the curve $\delta_k$ is equal to
  $T_{a_3}^k(\delta_0)$, where $a_3$ is as in Figure~\figurebasis.)
  The regular neighborhood of $\gamma\cup \delta_k$ is a lantern
  $L_k$, and we fix an identification of $L_k$ with our reference
  lantern $L$ by specifying that $\gamma$ and $\delta_k$ should
  correspond to $xy$ and $yz$ respectively. Let $f_k\in \Mag_g$ be the
  element corresponding under this identification to the mapping class
  $T_\alpha T_{\beta}^{-1}$ on $L$; it is easy to check using the
  lantern relation that $f_k$ is in fact
  $[T^{-1}_\gamma,T^{-1}_{\delta_k}]$. We will show that the images
  $\Psi(f_k)$ are linearly independent (over $\Z$).

  The boundary curves of $L_k$ are depicted in
  Figure~\figurefourcurves.
  \begin{figure}
    \label{figure:basis}
    \begin{center}
      \includegraphics[width=60mm]{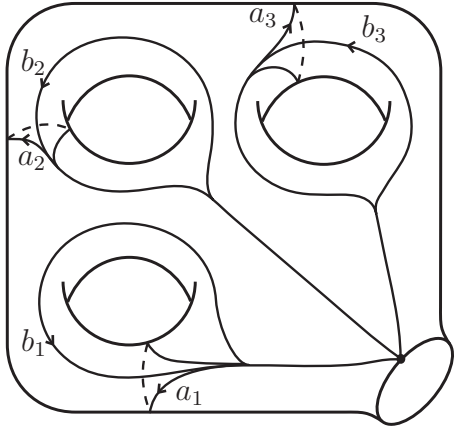}\\
    \end{center}
    \caption{A basis for $\pi_1(S_{g,1})$.}
  \end{figure}
  With the basis $a_1,b_1,\ldots,a_g,b_g$ for $\pi_1(S_{g,1})$ as illustrated in
  Figure~\figurebasis, we see that as curves $x$, $y$ and $z$
  can be represented by $[a_1,b_1]$, $[a_2,b_3a_3^kb_2]$, and
  $[b_2a_2^{-1}b_2^{-1}a_3,b_3a_3^k]$ respectively. As based loops, we
  actually have the conjugate
  $z={}^c[b_2a_2^{-1}b_2^{-1}a_3,b_3a_3^k]$, where
  $c=[b_3,a_3][b_2,a_2]a_2$. Note that with this representative for
  $z$, we have $xyz=[a_1,b_1][a_2,b_2][a_3,b_3]$, the fourth boundary
  curve in Figure~\figurefourcurves.
 
  Note that $a_1$ and $a_2$ intersect each $L_k$ in arcs corresponding
  to $A_1$ and $A_2$. Thus by Proposition~\ref{prop:lantern}, we have
  that \[\Psi(f_k)(\{a_1,a_2\})
  =(a_1-1)(a_2-1)\big[\big(\{a_1,b_1\}+\{a_2,b_3a_3^kb_2\}\big)\wedge
  a_2\{b_2a_2^{-1}b_2^{-1}a_3,b_3a_3^k\}\big]\] 
  Denote this element of $\bwedge^2 G^{\ab}$ by $\alpha_k$.  We now
  check that $\{\alpha_k\}$ is linearly independent as follows.  There
  is a standard embedding $G^{\ab}\hookrightarrow (\Z H)^{2g}$ given
  by sending the class $[x]$ to $(\partial x/\partial z_1,\ldots
  ,\partial x/\partial z_n)$, where $\{z_i\}$ is our basis for $F_n$
  and where $\partial/\partial z_i$ are the Fox derivatives (see e.g.\
  \cite{CP} for a detailed explanation of this embedding).  The only
  property of this embedding that we will need is that the components
  that make up $\alpha_k$ are mapped as follows by the embedding.
  Here the $A_i$ and $B_i$ make up a basis for $(\Z H)^{2g}$.
 \begin{align*}
   \{a_1,b_1\}\mapsto\quad&
   (1-b_1)A_1-(1-a_1)B_1\\
   \{a_2,b_3a_3^kb_2\}\mapsto\quad&(1-b_3a_3^kb_2)A_2
   \\-&(1-a_2)\big(B_3+b_3(1+\cdots+a_3^{k-1})A_3+b_3a_3^kB_2\big)\\
   \{b_2a_2^{-1}b_2^{-1}a_3,b_3a_3^k\}\mapsto\quad
   &(1-b_3a_3^k)\big((1-a_2^{-1})B_2-a_2^{-1}b_2A_2+a_2^{-1}A_3\big)\\
   -&(1-a_2^{-1}a_3)\big(B_3+b_3(1+\cdots+a_3^{k-1})A_3\big)
\end{align*}
By expanding out $\alpha_k$, we see that $\alpha_N$ is the only such
element which contains the term $A_1\wedge b_3a_3^NB_2$ with
nonzero coefficient; it follows that the $\alpha_k$ are linearly
independent, as desired.
\end{proof}

As the image of $\Psi$ is abelian, Theorem~\ref{theorem:imageA}
immediately implies Theorem~\ref{theorem:main} for $g\geq 3$.  Note
that the proof of Theorem~\ref{theorem:imageA} used in an essential
way that $g\geq 3$.  So in order to complete the proof of
Theorem~\ref{theorem:main}, we need another argument when $g=2$.

\begin{theorem}
  \label{theorem:imageB}
  $H_1(\Mag_2)$ has infinite rank.
\end{theorem}

\begin{proof} Suzuki showed that the element $f=[T_\gamma,T_\delta]$
  is in $\Mag_2$ for $\gamma$ and $\delta$ as in
  Figure~\figuregenustwo; in particular $\Mag_2$ is nontrivial. Let
  $S_2$ be a closed surface of genus 2; we denote by $\Tor_{2,\ast}$
  the Torelli group of $S_2$ with respect to a marked point $\ast$,
  and by $\Tor_2$ the Torelli group of the closed surface $S_2$. By
  Johnson \cite{J2}, we have the exact sequence \[1\to \Z\to
  \Tor_{2,1}\overset{p}{\to} \Tor_{2,\ast}\to 1,\] where the kernel is
  generated by a twist $T_\omega$ around the boundary $\omega=\partial
  S_2$. It is easy to check that the action of $T_\omega$ on
  $\pi_1(S_{2,1})$ is conjugation by $\omega$; since
  $\omega\not\in\Gamma^3$, we see that $T_\omega\not\in
  \Mag_2$. \begin{figure}
    \label{figure:genustwo}
    \begin{center}
      \includegraphics[width=120mm]{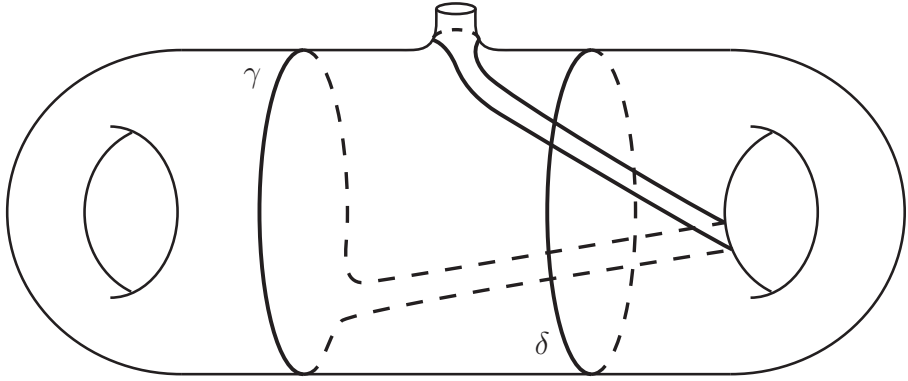}\\
    \end{center}
    \caption{The commutator $[T_\gamma,T_\delta]$ lies in $\Mag_2$.}
  \end{figure} It follows that $p$ restricts to an isomorphism between
  $\Mag_2$ and a subgroup $p(\Mag_2)<\Tor_{2,\ast}$.
  
  Again by Johnson \cite{J2}, we have the exact sequence
  \[1\to \Lambda \to \Tor_{2,\ast}\overset{\pi}{\to} \Tor_2\to 1,\]
  where $\Lambda\approx \pi_1(S_2,\ast)$; note that $\Tor_{2,\ast}$
  acts on $\pi_1(S_2,\ast)$, and the restriction to $\Lambda$ is just
  the action by conjugation. Mess \cite{Me} proved that $\Tor_2$ is
  free of infinite rank. It is easy to see from Figure~\figuregenustwo\ 
  that $f\in \ker \pi=\Lambda$. We use the following well-known lemma.
  \begin{lemma}\label{lemma}
    Any nontrivial infinite index normal subgroup of a surface group
    or free group is an infinite rank free group.
  \end{lemma}
  If $\pi\circ p(\Mag_2)<\Tor_2\approx F_\infty$ is nontrivial, then by
  Lemma~\ref{lemma}, $\Mag_2$ surjects to the infinite rank free group
  $\pi\circ p(\Mag_2)$, and we are done.

  Suppose that $p(\Mag_2)\subset \ker \pi=\Lambda$. Any $\varphi\in
  \Mag_2$ acts trivially on $\Gamma/\Gamma^3$; thus $p(\varphi)$ acts
  trivially on $\pi_1(S_2)/\pi_1(S_2)^3$. Since the action of
  $\Lambda$ is by conjugation, this implies that $p(\varphi)$ lies in
  $\Lambda^3$. Thus $p(\Mag_2)$ has infinite index in $\Lambda$, and
  so by Lemma~\ref{lemma}, $p(\Mag_2)\approx \Mag_2$ is an infinite
  rank free group.
\end{proof}

Theorem~\ref{theorem:main}, and hence Corollary~\ref{corollary:main},
follows immediately from Theorems~\ref{theorem:imageA} and
\ref{theorem:imageB}.

\begin{remark} One can check by explicit computation that for Suzuki's
  element $f\in \Mag_2$ above, $\Psi(f)=0$. It would be interesting to
  know whether $\Psi$ in fact vanishes on $\Mag_2$.
\end{remark}

\section{Computing the image of $\Phi$}
\label{section:braid}

The kernel $K$ of the map from $F_n=\langle x_1,\ldots,x_n\rangle$ to
$\Z=\langle t\rangle$ which sends each $x_i\mapsto t$ is normally
generated by the elements $x_ix_j^{-1}$. If we set $x_{i,k}\coloneq
x_1^kx_ix_1^{-k-1}$ for $i\neq 1$ and $k\in \Z$, then $\{x_{i,k}\}$
gives a basis for $K$ as a free group. As above, the conjugation of
$K$ by $F_n$ descends to a $\Z[t,t^{-1}]$ action on $K^{\ab}$. With
respect to this action we have $x_{i,k}=t^kx_{i,0}$, and thus
$K^{\ab}$ is a free $\Z[t,t^{-1}]$--module with basis
$\{y_i=x_{i,0}\}_{i\neq 1}$.

The braid group $B_n$ has generators $\sigma_1,\ldots,\sigma_{n-1}$;
the action of $\sigma_i$ on $F_n$ sends ${x_i\mapsto
  x_ix_{i+1}x_i^{-1}}$, ${x_{i+1}\mapsto x_i}$, and fixes the other
generators. The action of $B_n$ on $K^{\ab}$ commutes with the
$\Z[t,t^{-1}]$ action.

\begin{figure}[h]
  \label{figure:butterfly}
  \begin{center}
    \includegraphics[width=70mm]{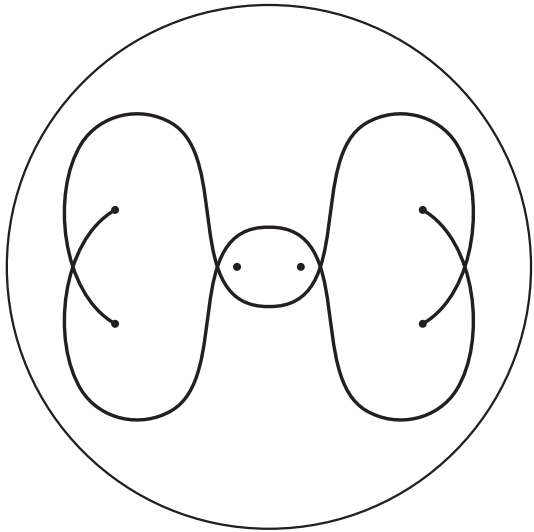}\\
  \end{center}
  \caption{The two arcs defining Bigelow's element $\phi_B$.}
\end{figure}

\begin{theorem}The image of $\Phi$ has infinite rank for $n\geq 6$.
\end{theorem}
\begin{proof}
  The element of the kernel found by Bigelow in \cite{Big} is the
  commutator of the half-twists along the arcs displayed in
  Figure~\figurebutterfly. In terms of the Artin generators,
  this is \[\phi_{B}=[\psi_1\sigma_3^{-1}\psi_1^{-1},\psi_2
  \sigma_3^{-1}\psi_2],\qquad\text{where } \psi_1=
  \sigma_4\sigma_5^{-1}\sigma_2^{-1}\sigma_1\quad\text{and }
  \psi_2=\sigma_4^{-1}\sigma_5^2\sigma_2\sigma_1^{-1}.\] In Appendix
  A, we give the computation of $\alpha\coloneq
  \Phi(\phi_{B})([x_2x_1^{-1}])=\Phi(\phi_{B})(y_2)$; it has 262
  terms. The only fact about $\alpha$ that we will need is that its highest term of
  the form $y_2\wedge t^ky_4$ is $-2y_2\wedge t^3y_4$, and its highest
  term of the form $y_2\wedge t^ky_5$ is $+2y_2\wedge t^2y_5$ (these
  terms are set in boxes in the appendix).

  It is easy to check that 
  \[\begin{array}{l}
    \sigma_4^2(x_4)=x_4x_5x_4x_5^{-1}x_4^{-1}\\
    \\
    \sigma_4^2(x_5)=x_5x_4x_5^{-1}\\
    \\
    \sigma_4^2(x_i)=x_i \ \ \text{for }i\neq 4,5.
  \end{array}\]
  By induction, for $k\geq 1$ we have
  \[\begin{array}{l}
    \sigma_4^{2k}(x_4)=(x_4x_5)^kx_4(x_4x_5)^{-k}\\
    \\
    \sigma_4^{2k}(x_5)=(x_4x_5)^{k-1}x_4x_5x_4^{-1}(x_4x_5)^{k-1}\\
    \\
    \sigma_4^{2k}(x_i)=x_i \ \ \text{for }i\neq 4,5.
  \end{array}\]
  The action of $\sigma_4^{2k}$ on $K^{\ab}$ in terms of our basis is
  thus given by:
  \begin{align*}
    y_4\ \ \mapsto&\ \ (1-t+t^2-\cdots-t^{k-1}+t^k)y_4\ \ +\ \
    (t-t^2+\cdots+t^{k-1}-t^k)y_5\\
    y_5\ \ \mapsto&\ \ (1-t+t^2-\cdots-t^{k-1})y_4\ \ \qquad+\ \
    (t-t^2+\cdots+t^{k-1})y_5\\
    y_i\ \ \mapsto&\ \ y_i\quad \text{for }i\neq 4,5
  \end{align*}
  Now for $k\geq 0$ set \[\alpha_k\coloneq
  \Phi(\sigma_4^{2k}\phi_{B}\sigma_4^{-2k})(y_2).\] By the
  equivariance of $\Phi$, and since $\sigma_4$ fixes $y_2$, we have
  $\alpha_k=\sigma_4^{2k}\cdot \alpha$. From the action of
  $\sigma_4^{2k}$ on $K^{\ab}$, we can see that the highest term in
  $\alpha_N$ of the form $y_2\wedge t^ky_4$ will be $-2y_2\wedge
  t^{3+N}y_4$.  Thus $\alpha_N$ is not contained in the span of
  $\{\alpha_1,\ldots,\alpha_{N-1}\}$; it follows that the $\alpha_k$
  are linearly independent over $\Z$, and thus the image of $\Phi$ has
  infinite rank.
\end{proof}
Theorem~\ref{theorem:main2} follows immediately.

\renewcommand{\thesection}{A}
\section{Appendix}
\label{section:appendix} The following computation was made, with the
method explained in Section \ref{section:braid}, with the help of 
\emph{Mathematica}.  A \emph{Mathematica} notebook implementing these
computations can be found at:

\bigskip
\url{http://math.uchicago.edu/~tchurch/infinitegeneration.html}
\bigskip

The output of this notebook is $\Phi(\phi_B)(y_2)$, which is:
\[\begin{array}{rrrrr}
-t^{-3}y_{2}\wedge t^{-2}y_{2}&+t^{-3}y_{2}\wedge t^{-1}y_{2}&-t^{-3}y_{2}\wedge y_{2}&-t^{-2}y_{2}\wedge y_{2}&+t^{-1}y_{2}\wedge y_{2}\\
+t^{-2}y_{2}\wedge ty_{2}&+t^{-1}y_{2}\wedge ty_{2}&-2y_{2}\wedge t^{2}y_{2}&+ty_{2}\wedge t^{3}y_{2}&+t^{2}y_{2}\wedge t^{3}y_{2}\\
-t^{3}y_{2}\wedge t^{4}y_{2}&+t^{-3}y_{2}\wedge t^{-4}y_{3}&-t^{-2}y_{2}\wedge t^{-4}y_{3}&-t^{-3}y_{2}\wedge t^{-3}y_{3}&+t^{-1}y_{2}\wedge t^{-3}y_{3}
\end{array}\]\[\begin{array}{rrrrr}
+t^{-2}y_{2}\wedge t^{-2}y_{3}&-t^{-1}y_{2}\wedge t^{-2}y_{3}&+t^{-3}y_{2}\wedge t^{-1}y_{3}&-y_{2}\wedge t^{-1}y_{3}&+ty_{2}\wedge t^{-1}y_{3}\\
-t^{2}y_{2}\wedge t^{-1}y_{3}&-2t^{-2}y_{2}\wedge y_{3}&+t^{3}y_{2}\wedge y_{3}&+t^{-1}y_{3}\wedge y_{3}&+2t^{-1}y_{2}\wedge ty_{3}\\
-t^{-1}y_{3}\wedge ty_{3}&-2y_{2}\wedge t^{2}y_{3}&-t^{4}y_{2}\wedge t^{2}y_{3}&+t^{-1}y_{3}\wedge t^{2}y_{3}&+ty_{2}\wedge t^{3}y_{3}\\
+t^{4}y_{2}\wedge t^{3}y_{3}&-y_{3}\wedge t^{3}y_{3}&+ty_{3}\wedge t^{3}y_{3}&-t^{2}y_{3}\wedge t^{3}y_{3}&+t^{-3}y_{2}\wedge t^{-3}y_{4}\\
-t^{-2}y_{2}\wedge t^{-3}y_{4}&-t^{-3}y_{2}\wedge t^{-2}y_{4}&+t^{-1}y_{2}\wedge t^{-2}y_{4}&+t^{-2}y_{2}\wedge t^{-1}y_{4}&-t^{-1}y_{2}\wedge t^{-1}y_{4}\\
+t^{-3}y_{2}\wedge y_{4}&-y_{2}\wedge y_{4}&+ty_{2}\wedge y_{4}&-t^{2}y_{2}\wedge y_{4}&-y_{3}\wedge y_{4}\\
+ty_{3}\wedge y_{4}&-t^{2}y_{3}\wedge y_{4}&-2t^{-2}y_{2}\wedge ty_{4}&+t^{3}y_{2}\wedge ty_{4}&+t^{-1}y_{3}\wedge ty_{4}\\
+t^{3}y_{3}\wedge ty_{4}&+y_{4}\wedge ty_{4}&+2t^{-1}y_{2}\wedge t^{2}y_{4}&-t^{-1}y_{3}\wedge t^{2}y_{4}&-t^{3}y_{3}\wedge t^{2}y_{4}\\
-y_{4}\wedge t^{2}y_{4}&\fbox{$-2y_{2}\wedge t^{3}y_{4}$}&-t^{4}y_{2}\wedge t^{3}y_{4}&+t^{-1}y_{3}\wedge t^{3}y_{4}&+t^{3}y_{3}\wedge t^{3}y_{4}\\
+y_{4}\wedge t^{3}y_{4}&+ty_{2}\wedge t^{4}y_{4}&+t^{4}y_{2}\wedge t^{4}y_{4}&-y_{3}\wedge t^{4}y_{4}&+ty_{3}\wedge t^{4}y_{4}\\
-t^{2}y_{3}\wedge t^{4}y_{4}&-ty_{4}\wedge t^{4}y_{4}&+t^{2}y_{4}\wedge t^{4}y_{4}&-t^{3}y_{4}\wedge t^{4}y_{4}&+t^{-3}y_{2}\wedge t^{-3}y_{5}\\
-t^{-2}y_{2}\wedge t^{-3}y_{5}&+t^{-3}y_{2}\wedge t^{-2}y_{5}&-t^{-2}y_{2}\wedge t^{-2}y_{5}&+y_{2}\wedge t^{-2}y_{5}&-t^{-4}y_{3}\wedge t^{-2}y_{5}\\
+t^{-3}y_{3}\wedge t^{-2}y_{5}&-t^{-1}y_{3}\wedge t^{-2}y_{5}&-t^{-3}y_{4}\wedge t^{-2}y_{5}&+t^{-2}y_{4}\wedge t^{-2}y_{5}&-y_{4}\wedge t^{-2}y_{5}\\
-t^{-3}y_{5}\wedge t^{-2}y_{5}&-2t^{-3}y_{2}\wedge t^{-1}y_{5}&+t^{-1}y_{2}\wedge t^{-1}y_{5}&+y_{2}\wedge t^{-1}y_{5}&-ty_{2}\wedge t^{-1}y_{5}\\
+t^{-4}y_{3}\wedge t^{-1}y_{5}&-t^{-2}y_{3}\wedge t^{-1}y_{5}&+2y_{3}\wedge t^{-1}y_{5}&+t^{-3}y_{4}\wedge t^{-1}y_{5}&-t^{-1}y_{4}\wedge t^{-1}y_{5}\\
+2ty_{4}\wedge t^{-1}y_{5}&+t^{-3}y_{5}\wedge t^{-1}y_{5}&+t^{-3}y_{2}\wedge y_{5}&+2t^{-2}y_{2}\wedge y_{5}&-2t^{-1}y_{2}\wedge y_{5}\\
-y_{2}\wedge y_{5}&-t^{2}y_{2}\wedge y_{5}&-t^{-3}y_{3}\wedge y_{5}&+t^{-2}y_{3}\wedge y_{5}&-y_{3}\wedge y_{5}\\
-ty_{3}\wedge y_{5}&-t^{2}y_{3}\wedge y_{5}&-t^{-2}y_{4}\wedge y_{5}&+t^{-1}y_{4}\wedge y_{5}&-ty_{4}\wedge y_{5}\\
-t^{2}y_{4}\wedge y_{5}&-t^{3}y_{4}\wedge y_{5}&+t^{-1}y_{5}\wedge y_{5}&-t^{-3}y_{2}\wedge ty_{5}&-t^{-1}y_{2}\wedge ty_{5}\\
+y_{2}\wedge ty_{5}&+ty_{2}\wedge ty_{5}&+t^{3}y_{2}\wedge ty_{5}&+t^{-1}y_{3}\wedge ty_{5}&-y_{3}\wedge ty_{5}\\
+2ty_{3}\wedge ty_{5}&+t^{3}y_{3}\wedge ty_{5}&+y_{4}\wedge ty_{5}&-ty_{4}\wedge ty_{5}&+2t^{2}y_{4}\wedge ty_{5}\\
+t^{4}y_{4}\wedge ty_{5}&-t^{-2}y_{5}\wedge ty_{5}&-y_{5}\wedge ty_{5}&+t^{-2}y_{2}\wedge t^{2}y_{5}&-t^{-1}y_{2}\wedge t^{2}y_{5}\\
\fbox{$+2y_{2}\wedge t^{2}y_{5}$}&-t^{2}y_{2}\wedge t^{2}y_{5}&+t^{3}y_{2}\wedge t^{2}y_{5}&-t^{-1}y_{3}\wedge t^{2}y_{5}&-t^{2}y_{3}\wedge t^{2}y_{5}\\
-y_{4}\wedge t^{2}y_{5}&-t^{3}y_{4}\wedge t^{2}y_{5}&+t^{-1}y_{5}\wedge t^{2}y_{5}&-2y_{5}\wedge t^{2}y_{5}&+ty_{5}\wedge t^{2}y_{5}\\
-ty_{2}\wedge t^{3}y_{5}&-t^{2}y_{2}\wedge t^{3}y_{5}&-t^{4}y_{2}\wedge t^{3}y_{5}&+y_{3}\wedge t^{3}y_{5}&+ty_{4}\wedge t^{3}y_{5}\\
+ty_{5}\wedge t^{3}y_{5}&+t^{2}y_{5}\wedge t^{3}y_{5}&+t^{2}y_{2}\wedge t^{4}y_{5}&+t^{3}y_{2}\wedge t^{4}y_{5}&-ty_{3}\wedge t^{4}y_{5}\\
+t^{3}y_{3}\wedge t^{4}y_{5}&-t^{2}y_{4}\wedge t^{4}y_{5}&+t^{4}y_{4}\wedge t^{4}y_{5}&-t^{2}y_{5}\wedge t^{4}y_{5}&-t^{3}y_{5}\wedge t^{4}y_{5}\\
-t^{3}y_{2}\wedge t^{5}y_{5}&+t^{2}y_{3}\wedge t^{5}y_{5}&-t^{3}y_{3}\wedge t^{5}y_{5}&+t^{3}y_{4}\wedge t^{5}y_{5}&-t^{4}y_{4}\wedge t^{5}y_{5}\\
+t^{3}y_{5}\wedge t^{5}y_{5}&-t^{-3}y_{2}\wedge t^{-3}y_{6}&+t^{-2}y_{2}\wedge t^{-3}y_{6}&-t^{-2}y_{5}\wedge t^{-3}y_{6}&+t^{-1}y_{5}\wedge t^{-3}y_{6}\\
+t^{-3}y_{2}\wedge t^{-2}y_{6}&-t^{-1}y_{2}\wedge t^{-2}y_{6}&+t^{-2}y_{5}\wedge t^{-2}y_{6}&-y_{5}\wedge t^{-2}y_{6}&+t^{-3}y_{2}\wedge t^{-1}y_{6}\\
-t^{-2}y_{2}\wedge t^{-1}y_{6}&+y_{2}\wedge t^{-1}y_{6}&-t^{-4}y_{3}\wedge t^{-1}y_{6}&+t^{-3}y_{3}\wedge t^{-1}y_{6}&-t^{-1}y_{3}\wedge t^{-1}y_{6}\\
-t^{-3}y_{4}\wedge t^{-1}y_{6}&+t^{-2}y_{4}\wedge t^{-1}y_{6}&-y_{4}\wedge t^{-1}y_{6}&-t^{-3}y_{5}\wedge t^{-1}y_{6}&+ty_{5}\wedge t^{-1}y_{6}\\
+t^{-3}y_{6}\wedge t^{-1}y_{6}&-t^{-2}y_{6}\wedge t^{-1}y_{6}&-t^{-3}y_{2}\wedge y_{6}&-t^{-2}y_{2}\wedge y_{6}&+t^{-1}y_{2}\wedge y_{6}\\
+2y_{2}\wedge y_{6}&-2ty_{2}\wedge y_{6}&+t^{2}y_{2}\wedge y_{6}&+t^{-3}y_{3}\wedge y_{6}&-t^{-2}y_{3}\wedge y_{6}\\
-t^{-1}y_{3}\wedge y_{6}&+3y_{3}\wedge y_{6}&-ty_{3}\wedge y_{6}&+t^{2}y_{3}\wedge y_{6}&+t^{-2}y_{4}\wedge y_{6}\\
-t^{-1}y_{4}\wedge y_{6}&-y_{4}\wedge y_{6}&+3ty_{4}\wedge y_{6}&-t^{2}y_{4}\wedge y_{6}&+t^{3}y_{4}\wedge y_{6}\\
-y_{5}\wedge y_{6}&+ty_{5}\wedge y_{6}&-2t^{2}y_{5}\wedge y_{6}&-t^{-2}y_{6}\wedge y_{6}&+t^{-3}y_{2}\wedge ty_{6}\\
+t^{-2}y_{2}\wedge ty_{6}&-y_{2}\wedge ty_{6}&-ty_{2}\wedge ty_{6}&-t^{3}y_{2}\wedge ty_{6}&-t^{-1}y_{3}\wedge ty_{6}\\
+y_{3}\wedge ty_{6}&-2ty_{3}\wedge ty_{6}&-t^{3}y_{3}\wedge ty_{6}&-y_{4}\wedge ty_{6}&+ty_{4}\wedge ty_{6}\\
-2t^{2}y_{4}\wedge ty_{6}&-t^{4}y_{4}\wedge ty_{6}&+t^{-2}y_{5}\wedge ty_{6}&+t^{-1}y_{5}\wedge ty_{6}&+t^{2}y_{5}\wedge ty_{6}\\
+t^{3}y_{5}\wedge ty_{6}&+t^{-1}y_{6}\wedge ty_{6}&+2y_{6}\wedge ty_{6}&-t^{-2}y_{2}\wedge t^{2}y_{6}&-t^{-1}y_{2}\wedge t^{2}y_{6}\\
+t^{2}y_{2}\wedge t^{2}y_{6}&+t^{-1}y_{3}\wedge t^{2}y_{6}&+t^{2}y_{3}\wedge t^{2}y_{6}&+t^{3}y_{3}\wedge t^{2}y_{6}&+y_{4}\wedge t^{2}y_{6}\\
+t^{3}y_{4}\wedge t^{2}y_{6}&+t^{4}y_{4}\wedge t^{2}y_{6}&-t^{-1}y_{5}\wedge t^{2}y_{6}&+ty_{5}\wedge t^{2}y_{6}&-t^{2}y_{5}\wedge t^{2}y_{6}\\
-t^{4}y_{5}\wedge t^{2}y_{6}&-2y_{6}\wedge t^{2}y_{6}&-ty_{6}\wedge t^{2}y_{6}&+2y_{2}\wedge t^{3}y_{6}&+t^{4}y_{2}\wedge t^{3}y_{6}
\end{array}\]\[\begin{array}{rrrrr}
-t^{-1}y_{3}\wedge t^{3}y_{6}&-t^{3}y_{3}\wedge t^{3}y_{6}&-y_{4}\wedge t^{3}y_{6}&-t^{4}y_{4}\wedge t^{3}y_{6}&-y_{5}\wedge t^{3}y_{6}\\
-t^{2}y_{5}\wedge t^{3}y_{6}&+t^{5}y_{5}\wedge t^{3}y_{6}&+y_{6}\wedge t^{3}y_{6}&+t^{2}y_{6}\wedge t^{3}y_{6}&-ty_{2}\wedge t^{4}y_{6}\\
-t^{2}y_{2}\wedge t^{4}y_{6}&-t^{4}y_{2}\wedge t^{4}y_{6}&+y_{3}\wedge t^{4}y_{6}&+ty_{4}\wedge t^{4}y_{6}&+ty_{5}\wedge t^{4}y_{6}\\
+t^{2}y_{5}\wedge t^{4}y_{6}&+t^{4}y_{5}\wedge t^{4}y_{6}&-t^{5}y_{5}\wedge t^{4}y_{6}&-ty_{6}\wedge t^{4}y_{6}&+t^{3}y_{2}\wedge t^{5}y_{6}\\
-t^{2}y_{3}\wedge t^{5}y_{6}&+t^{3}y_{3}\wedge t^{5}y_{6}&-t^{3}y_{4}\wedge t^{5}y_{6}&+t^{4}y_{4}\wedge t^{5}y_{6}&-t^{3}y_{5}\wedge t^{5}y_{6}\\
+t^{3}y_{6}\wedge t^{5}y_{6}&-t^{4}y_{6}\wedge t^{5}y_{6}
\end{array}\]

\small

\noindent
Dept. of Mathematics, University of Chicago\\
5734 University Ave.\\
Chicago, IL 60637\\
E-mails: tchurch@math.uchicago.edu, farb@math.uchicago.edu

\end{document}